\documentclass[a4paper,11pt]{article}
\usepackage{graphicx}
\usepackage{amssymb}
\usepackage{latexsym, bm}
\usepackage{multicol}
\usepackage{indentfirst}
\usepackage{amssymb,amsfonts}
\usepackage{amsmath}
\usepackage{setspace}
\usepackage{amssymb}
\newtheorem{theorem}{Theorem}[section]
\newtheorem{lemma}{Lemma}[section]
\newtheorem{claim}{Claim}[section]

\newtheorem{corollary}{Corollary}[section]

\newtheorem{definition}{Definition}[section]

\newtheorem{conjecture}{Conjecture}[section]

\begin{spacing}{1}
\begin{document}

\title{Ramsey non-goodness involving books}
\date{}

\author{
Chunchao Fan\footnote{Center for Discrete Mathematics, Fuzhou University, Fuzhou, 350108 P.~R.~China. Email: {\tt 1807951575@qq.com}.} \;\; and\;\;
Qizhong Lin\footnote{Corresponding author. Center for Discrete Mathematics, Fuzhou University,
Fuzhou, 350108 P.~R.~China. Email: {\tt linqizhong@fzu.edu.cn}. Supported in part by NSFC(No.\ 1217010182).}
}

\maketitle
\begin{abstract}

In 1983, Burr and Erd\H{o}s  initiated the study of Ramsey goodness problems.
Nikiforov and Rousseau (2009) resolved almost all goodness questions raised by Burr and Erd\H{o}s, in which the bounds on the parameters are of tower type since their proofs rely on the regularity lemma.
Let $B_{k,n}$ be the book graph on $n$ vertices which consists of $n-k$ copies of $K_{k+1}$ all sharing a common $K_k$,  and let $H=K_p(a_1,\dots,a_{p})$ be the complete $p$-partite graph with parts of sizes $a_1,\dots,a_{p}$.

Recently, avoiding use of the regularity lemma, Fox, He and Wigderson (2021) revisit several Ramsey goodness results involving books.
They comment that it would be very interesting to see how far one can push these ideas.
In particular, they conjecture that for all integers $k, p, t\ge 2$, there exists some $\delta>0$ such that for all $n\ge 1$, $1\leq a_1\le\cdots\le a_{p-1}\le t$ and $a_p \le \delta n$, we have $r(H, B_{k,n})= (p-1)(n-1)+d_k(n,K_{a_1,a_2})+1,$ where $d_k(n,K_{a_1,a_2})$ is the maximum $d$ for which there is an $(n+d-1)$-vertex $K_{a_1,a_2}$-free graph in which at most $k-1$ vertices have degree less than $d$.
They verify the conjecture when $a_1=a_2=1$.

Building upon the work of Fox et al. (2021), we make a substantial step by showing that the conjecture  ``roughly" holds if $a_1=1$ and $a_2|(n-1-k)$, i.e. $a_2$ divides $n-1-k$.
Moreover, we prove that for every $k, a\geq 1$ and $p\ge2$, there exists $\delta>0$ such that for all  large $n$ and $b\le  \delta\ln n$,  $r(K_p(1,a,b,\dots,b), B_{k,n})= (p-1)(n-1)+k(p-1)(a-1)+1$ if $a|(n-1-k)$, where the case when $a=1$ has been proved by Nikiforov and Rousseau (2009) using the regularity lemma.  The bounds on $1/\delta$ we obtain are not of tower-type since our proofs do not rely on the regularity lemma.

\medskip

\end{abstract}

\section{Introduction}

For graphs $H_1$ and $H_2$, the Ramsey number $r(H_1,H_2)$ is the smallest positive integer $N$ such that any graph $G$ on $N$ vertices contains $H_1$ as a subgraph, or its complement contains $H_2$ as a subgraph. The existence of the Ramsey number $r(H_1,H_2)$ follows from \cite{ram}.

The following  result of Chv\'atal is well-known.
\begin{theorem}[Chv\'atal \cite{c}] \label{chv}
For any tree $T_n$ of order $n$,  $r(K_p,T_n)= (p-1)(n-1)+1.$
\end{theorem}

For a graph $H$, let  $\chi(H)$ be the chromatic number of $H$, and let $s(H)$ be the minimum size of a color class over all proper vertex-colorings of $H$ by $\chi(H)$ colors. Burr \cite{bur} first observed that if $H_1$ is a graph with chromatic number $\chi(H_1)=p$  and  $H_2$ is a connected graph of order $n\geq s(H_1)$, then
\begin{align}\label{burr}
r(H_1,H_2)\ge(p-1)(n-1)+s(H_1).
\end{align}

In 1983, Burr and Erd\H{o}s \cite{be} initiated the study of Ramsey goodness problems.  We say that a connected graph $H_2$ is {\em $H_1$-good} if the equality of (\ref{burr}) holds.
Thus all trees are $K_p$-good from the above result by Chv\'{a}tal \cite{c}.
Ramsey goodness problems have attracted a great deal of attention,  we refer the reader to the survey \cite[Section 2.5]{cfs-15} and many other related references.

Let $B_{k,n}$ (or $B_{n-k}^{(k)}$) denote the book graph on $n$ vertices which consists of $n-k$ copies of $K_{k+1}$ all sharing a common $K_k$,  and let $K_p(a_1,\dots,a_{p})$ be the complete $p$-partite graph with parts of sizes $a_1,\dots,a_{p}$.  For convenience, we denote by $B_n$ and $K_{m,n}$ instead of $B_{2,n+2}$ and $K_2(m,n)$, respectively. Also, we will write $K_p(a)$ for $K_p(a,\dots,a)$.

Using the regularity lemma \cite{sz}, Nikiforov and Rousseau \cite{nr} showed that large books are $K_{p}$-good.
Subsequently, extending the method used in \cite{nr}, Nikiforov and Rousseau  \cite{nr09} obtained extremely general goodness results, in which they resolved almost all goodness questions raised by Burr and Erd\H{o}s. As a special case, they proved the following result.
\begin{theorem}[Nikiforov and Rousseau  \cite{nr09}]\label{nr-b}
For all $k,p\ge2$, there exists some $\delta>0$ such that  for all large $n$, $$r(B_{p-1,\delta n},B_{k,n})=(p-1)(n-1)+1.$$
\end{theorem}

All the bounds on $n$ of these results are of tower type since the proofs rely on the regularity lemma, and Nikiforov and Rousseau raised the question of what the best possible $n$ is.
Recently, avoiding use of the regularity lemma, Fox, He and Wigderson \cite{fox} obtain that if $n\ge2^{k^{10p}}$, then $B_{k,n}$ is already $K_p$-good. Moreover, they establish the following result. 
\begin{theorem}[Fox, He and Wigderson \cite{fox}]\label{fox}
For every $k, p, t\geq 2$, there exists $\delta>0$ such that the following holds for all large $n$. Let $1\leq a_1\leq \dots \leq a_{p-1}\leq t$ and $a_{p}\leq\delta n$ be positive integers. If $a_1=a_2=1$, then $r(K_p(a_1,\dots,a_{p}), B_{k,n})= (p-1)(n-1)+1$.
\end{theorem}

The main novelty is that the proof again avoids the use of the regularity lemma, so that the bounds on $1/\delta$ are not of tower-type.
Indeed, once $n$ is double-exponential in $k, p$, and $t$, the statement holds with $1/\delta$ merely single-exponential in $k, p$, and $t$.

Based on the above result, they ask to determine $r(K_p(a_1,\dots,a_{p}), B_{k,n})$ for larger values of $a_1$ and $a_2$.
They also comment that it would be very interesting to see how far one can push these ideas.
The following is the so-called Dirac-type extremal function.
\begin{definition}[Fox, He and Wigderson \cite{fox}]\label{def-fhw}
Given a graph $H$ and integers $k, n$, let $d_k(n, H)$ be the maximum $d$ for which
there is an $(n + d-1)$-vertex $H$-free graph, at most $k-1$ vertices of which have degree less than $d$.
\end{definition}

Let $d=d_k(n,K_{a_1,a_2})$, and let $G_1$ be a $K_{a_1,a_2}$-free graph on $n+d-1$ vertices, at most $k-1$ of which have degree less than $d$.
It follows that at most $k-1$ vertices in $\overline{G_1}$ have degree at least $n-1$, and so $\overline{G_1}$ contains no $B_{k,n}$. Let $G$ be a graph on $(p-1)(n-1)+d$ vertices, whose vertex set $V(G)=\sqcup_{1\le \i\le p-1}U_i$ with $|U_1|=n+d-1$ and $|U_2|=\cdots=|U_{p-1}|=n-1$, such that $G[U_1]$ is isomorphic to $G_1$, and all pairs of vertices in different parts are adjacent, and $U_i$ forms an independent set for $2\le i\le p-1$. It is clear that $G$ contains no $K_p(a_1,\dots,a_{p})$ and its complement $\overline{G}$ contains no $B_{k,n}$. Therefore,
\begin{align}\label{obs-f}
r(K_p(a_1,\dots,a_{p}), B_{k,n})>(p-1)(n-1)+d_k(n,K_{a_1,a_2}).
\end{align}

From the above observation, they propose the following conjecture.

\begin{conjecture}[Fox, He and Wigderson \cite{fox}]\label{conj}
For every $k, p, t\geq 2$, there exists $\delta>0$ such that the following holds for all $n\ge1$. For $1\le a_1\leq \dots \leq a_{p-1}\leq t$ and $a_{p}\leq\delta n$,
$$r(K_p(a_1,\dots,a_{p}), B_{k,n})=(p-1)(n-1)+d_k(n,K_{a_1, a_2})+1.$$
\end{conjecture}

Fox, He and Wigderson verify Conjecture \ref{conj} for $a_1=a_2=1$, see \cite[Theorem 1.3]{fox} or the above Theorem \ref{fox}.
Note that $d_k(n,K_{a_1, a_2})=0$ if $a_1=a_2=1$. 

In this paper, we first consider Conjecture \ref{conj} for $a_1=1$ and $a_2\ge1$.  Building upon the work of Fox, He and Wigderson \cite{fox}, we show the following result avoiding use of the regularity lemma.
For positive integers $a$ and $b$, we write $a|b$ if $a$ divides $b$.
\begin{theorem}\label{main}
For every $k, p, t\geq 2$, there exists $\delta>0$ such that the following holds for all large $n$. Let $1\leq a_1\leq \dots \leq a_{p-1}\leq t$ and $a_{p}\leq\delta n$ be positive integers. If $a_1=1$, then $$r(K_p(a_1,\dots,a_{p}), B_{k,n})\le (p-1)(n-1)+k(p-1)(a_2-1)+1.$$
The equality holds if $a_2|(n-1-k)$.
\end{theorem}

Theorem \ref{main} implies that Conjecture \ref{conj} ``roughly'' holds when $a_1=1$ and $a_2|(n-1-k)$. However, we shall see that Conjecture \ref{conj} does not hold for many cases in the next section. Moreover, the above result strengthens \cite[Theorem 1]{lin} in which the multipartite graph is fixed and the regularity lemma is a key tool.

\medskip
Moreover, avoiding use of the regularity lemma, we obtain the following result.

\begin{theorem}\label{main-2}
For every $k, a\geq 1$ and $p\ge2$, there exists $\delta>0$ such that for all large $n$ and $b\le\delta\ln n$,
$r(K_p(1,a,b,\dots,b), B_{k,n})\le (p-1)(n-1)+k(p-1)(a-1)+1.$
The equality holds if $a|(n-1-k)$.
\end{theorem}

Let us point out that since our proofs do not rely the regularity lemma, all the bounds on $1/\delta$ are not of tower-type.
Also, we may take $n$ as double-exponential in $k, p$, and $t$ (or $a$).

\medskip
As a corollary, we have the following result which has been obtained by  Nikiforov and Rousseau  \cite[Theorem 2.2]{nr09}. Indeed, they obtained a much more general result in which the book graph can be replaced by a family of sparse graphs.
\begin{corollary}[Nikiforov and Rousseau  \cite{nr09}]\label{nr-lnn}
For every $k\geq 1$ and $p\ge2$, there exists $\delta>0$ such that
$r(K_p(1,1,\lceil \delta\ln n\rceil,\dots,\lceil \delta\ln n\rceil), B_{k,n})= (p-1)(n-1)+1$
 for all large $n$.
\end{corollary}


\medskip
The rest of the paper is organized as follows. In Section \ref{disp}, we will see that Conjecture \ref{conj} does not hold in many cases. In Section \ref{pre}, we will collect several stability result including a stability-supersaturation lemma, and we will also obtain stable structures with forbidden graphs. In Section \ref{pf-1} and Section \ref{pf-2}, we shall present the proofs of Theorem \ref{main} and Theorem \ref{main-2}.  Finally, we will mention some interesting problems in Section \ref{clu}.

\section{Conjecture \ref{conj} does not hold in many cases}\label{disp}


 In this section, let us have an overview on Conjecture \ref{conj} proposed by Fox, He and Wigderson \cite{fox}.
It would be very difficult to draw the whole picture of this conjecture, since it contains many special cases that are notorious difficult themselves.

  \medskip

1. In the following, we will see that  Conjecture \ref{conj} does not hold in many cases.
 For $a_1=1$, let $G$ be a $K_{1,a_2}$-free graph on $n+d-1$ vertices, at most $k-1$ of which have degree less than $d$. It follows that $G$ contains at least $n+d-k$ vertices with degree at least $d$, implying $d\leq a_2-1$. From Definition \ref{def-fhw}, we have $d_k(n,K_{1, a_2})\le a_2-1$. Therefore, if Conjecture \ref{conj} holds, then we must have an upper bound as
 $$r(K_p(a_1,\dots,a_{p}), B_{k,n})\leq (p-1)(n-1)+(a_2-1)+1,$$
 which will lead to a contradiction from the following construction.

Indeed, for $a_1=1$ and $a_2,k, p\geq 2$, 
 let $F$ be $\lfloor (n-k-1)/a_2\rfloor+k$ disjoint copies of $K_{a_2}$, and let $G$ be the graph obtained from $p-1$ disjoint copies of $F$ by adding all edges between them. We can see that $G$ contains no $K_p(1,a_2,\dots,a_{p})$ for $1=a_1<a_2\le \dots\le a_p$, and its complement $\overline{G}$ contains no $B_{k,n}$ since the number of the common neighbors of any clique of $k$ vertices in $\overline{G}$ is
 $(\lfloor \frac {n-k-1}{a_2}\rfloor+k)\cdot a_2-ka_2\le n-k-1.$ Therefore,
\begin{align}\label{low-bou}
r(K_p(1,a_2,\dots,a_{p}), B_{k,n})&>(p-1)\left(\left\lfloor \frac {n-k-1}{a_2}\right\rfloor+k\right)a_2
\\&\ge (p-1)\left( \frac {n-k-1}{a_2}-\frac {a_2-1}{a_2}+k\right)a_2\nonumber
\\&=(p-1)(n-1)+(p-1)(k-1)(a_2-1),\nonumber
\end{align}
which is at least $(p-1)(n-1)+(a_2-1)+1$ provided $a_2,k, p\geq 2$ and $(k,p)\neq(2,2)$.
 Consequently, Conjecture \ref{conj} does not hold for every $k, p\geq 2$ with $(k,p)\neq(2,2)$, $a_1=1$ and $2\leq a_2\leq \dots \leq a_{p}$.
In particular, the above construction implies that the book graph $B_{k,n}$ is not $K_p(1,a_2,\dots,a_{p})$-good for $a_2\ge2$.

\medskip
2. When $a_1\ge2$, let us see that many special cases of Conjecture \ref{conj} hold when $p=2$. The case when $p=2$ may be referred to as the ``degenerate'' case of goodness problems.
 Recall the well-known Erd\H{o}s-Renyi graph $ER_q$: for any prime  power $q$, there exists a $C_4$-free graph of order $q^2+q+1$ whose minimum degree is $q$,  see Erd\H{o}s, Renyi and S\'{o}s \cite{ER66} (or Brown \cite{Br66}). Therefore, we have
 $r(C_4,K_{1,\,q^2+1})> q^2+q+1$. Since $ER_q$ has no vertex of degree less than $q$, from Definition \ref{def-fhw}, we have $d_1(q^2+2,C_4)\ge q$. Parsons \cite{par} showed that $r(C_4,K_{1,\,q^2+1})= q^2+q+2$ for any prime  power $q$, implying $d_1(q^2+2,C_4)\le q$.
Therefore, Conjecture \ref{conj} is true for this special case. For $k=1$, one can also see \cite{par2,zcc} that many other cases of Conjecture \ref{conj} hold.

For the special case $r(C_4,B_{n})$, we only know the exact value of $r(C_4,B_{n})$ for $1\le n\le 14$ (see \cite{cs,frs,ts1,ts2}) until recently the second author with Li and Peng \cite{llp} obtain the exact value of $r(C_4,B_{n})$ for infinitely many $n$.
In particular, we show $r(C_4,B_{q^2-q-2}) = q^2+q-1$ for any prime power $q\ge4$, in which the lower bound follows since there exists a $C_4$-free graph $G_q$ of order $q^2+q-2$ whose minimum degree is $q$. However, such graph is not what we want from Conjecture \ref{conj}. Indeed, since $q^2+q-2=(q^2-q-1)+(2q-1)$, by Definition \ref{def-fhw}, we need to find a $C_4$-free graph on $q^2+q-2$ vertices in which at most one vertex has degree less than $2q-1$, which does not satisfied by $G_q$.
Conjecture \ref{conj} is still open for other cases when $a_1\ge2$. Clearly, all large books are not $C_4$-good from the lower bound $r(C_4,K_{1,n})\ge n+\lfloor n^{1/2}-6n^{11/40}\rfloor$ by Burr et al. \cite{befrs}.

\section{Preliminaries}\label{pre}

Let $H=(V,E)$ be a graph. For any disjoint subsets $A,B\subseteq V$, we use $e(A,B)$ to denote the number of edges between $A$ and $B$.
Also, we use $e(A)$ to denote the number of edges in $A$.
Let $H[A]$ be the subgraph induced by $A$.
For a vertex $v\in V$, let $d_H(v)$ denote its degree in graph $H$ and let $\delta(H)$ denote its minimum degree among all vertices in $H$.
We denote by $N_A(v)$ the neighborhood of $v$ in $A$, and let $d_A(v)=|N_A(v)|$.

\subsection{Stability results}
The following result by Tur\'{a}n \cite{turan}  guarantees a large independent set in a sparse graph.
\begin{lemma} [Tur\'{a}n \cite{turan}]\label{turan}
For any graph $G$ of order $n$ with average degree $d$, the independence number
$\alpha(G)$ is at least $\frac{n}{1+d}.$
\end{lemma}

The following lemma is essentially due to Erd\H{o}s \cite{erdos64}, which implies that if a graph $G$ with $n$ vertices and $\Omega(n^p)$ copies of $K_p$, then we can find a copy of $K_p(a_1,\dots,a_{p})$ with one part of size linear on $n$. One can also see Fox, He and Wigderson \cite{fox} for a proof.
\begin{lemma}[Erd\H{o}s \cite{erdos64}]\label{Erdos}
For every $\eta>0$, and $p,t\geq2$, and $1\leq a_1\leq \dots \leq a_{p-1}\leq t$, there exists some $\delta>0$ such that the following holds for large enough $n$. If $G$ is a $K_p(a_1,\dots,a_{p})$-free graph on $n$ vertices with $a_p\leq \delta n$, then $G$ has at most $\eta n^p$ copies of $K_p$.
\end{lemma}

The well-known Erd\H{o}s-Simonovits stability theorem \cite{er66,sim} states that if a $K_p$-free graph has slightly fewer edges than the $K_p$-free Tur\'{a}n graph, then its structure is very close to the Tur\'{a}n graph. Moreover, a supersaturation result of  Erd\H{o}s-Simonovits \cite{es} says that if a $K_p$-free graph on $n$ vertices has slightly more edges than the $K_p$-free Tur\'{a}n graph, then it must contain $\Omega(n^p)$ copies of  $K_p$.

The following {\em stability-supersaturation lemma} by Fox, He and Wigderson \cite[Theorem 3.1]{fox} (in a slightly different form) is one of the key ingredients which will be used to prove the desired structure for graphs forbidding some special classes of graphs, which can be referred to as a combination of the stability theorem \cite{er66,sim} and  the supersaturation result \cite{es}. This stability-supersaturation lemma implies that if a graph $G$ has slightly fewer edges than the $K_p$-free Tur\'{a}n graph and has few copies of $K_p$, then it is close to the Tur\'{a}n graph.

\begin{lemma}[Fox, He and Wigderson \cite{fox}]\label{stability}
For every $\varepsilon>0$ and every integer $p\geq 2$, there exist $\eta,\gamma>0$ such that the following holds for all $n\geq 10$. Suppose $G$ is a graph on $n$ vertices with minimum degree at least $(1-\frac{1}{p-1}-\gamma)n$ and at most $\eta n^p$ copies of $K_p$. Then there is a partition $V(G)=\sqcup_{1\le i\le p-1} V_i$, such that the following hold:

\medskip
(i) $\sum_{1\le i\le p-1}e(V_i)\le\varepsilon {n\choose 2}$.

(ii) $\big||V_i|-\frac{n}{p-1}\big|\leq\sqrt{2\varepsilon}n$.

(iii) $e(V_i,V_j)\geq (1-p^2\varepsilon)|V_i||V_j|$.

(iv) For each $v\in V_i$, $d_{V_i}(v)\leq d_{V_j}(v)$.

\medskip
Moreover, we may take $\gamma=\min\{\frac{1}{2p^2},\frac{\varepsilon}{2}\}$ and $\eta=p^{-10p}\varepsilon$.
\end{lemma}

One can see a stronger result (weakening the minimum degree condition to an average degree condition) due to Conlon, Fox and Sudakov
\cite[Corollary 3.4]{cfs} by using the celebrated graph removal lemma (see e.g. Conlon and Fox \cite{cf}), however, which will give tower-type bounds in the parameters.

\subsection{Stable structures with forbidden graphs}
In \cite{aes}, Andr\'{a}sfai-Erd\H{o}s-S\'{o}s established a minimum-degree stability version of Tur\'{a}n's theorem, which states that if an $n$-vertex $K_p$-free graph has minimum degree greater than $\frac{3p-7}{3p-4}n$, then it is $(p-1)$-partite. Here the constant $\frac{3p-7}{3p-4}$ is best possible. Extending a result of Fox, He and Wigderson \cite[Theorem 4.1]{fox}, we obtain the following lemma,  which states that if a $K_p(1,a_2,\dots,a_{p})$-free graph has high minimum degree, then it is nearly $(p-1)$-partite, in which the number of internal edges in each part is bounded. 
\begin{lemma}\label{AES Th}
For every $p, t\geq2$, every $1\leq a_2\leq \dots \leq a_{p-1}\leq t$, there exist some $\gamma,\delta>0$ such that if $n$ is large enough in terms of $t$ and $p$, $a_p\leq \delta n$, and $G$ is a $K_p(1,a_2,\dots,a_{p})$-free graph on $n$ vertices with minimum degree at least $(1-\frac{1}{p-1}-\gamma)n$, then there is a partition $V(G)=\sqcup_{1\le i\le p-1} V_i$ such that each vertex of $V_i$ has at most $a_2-1$ neighbors in $V_i$, for $i=1, 2, \dots, p-1$.
\end{lemma}
\noindent{\bf Proof.} 
 For $\varepsilon>0$, where $\varepsilon$ is sufficiently small in terms of $p$ and $t$, let $\gamma$ and $\eta$ be applied in Lemma \ref{stability}. Further, let $\delta>0$ be the parameter guaranteed by Lemma \ref{Erdos}. For $K_p(1,a_2,\dots,a_{p})$-free graph $G$, we have that $G$ has at most $\eta n^p$ copies of $K_p$ by Lemma \ref{Erdos}. Since $G$ has at most $\eta n^p$ copies of $K_p$ with minimum degree at least $(1-\frac{1}{p-1}-\gamma)n$, it follows from  Lemma \ref{stability} that there is a partition $V(G)=\sqcup_{1\le i\le p-1} V_i$ such that (i)-(iv) of Lemma \ref{stability} hold.

\begin{claim}\label{ext-fhw}
For distinct $i$ and $j$, $1\leq i,j\leq p-1$, we have

\medskip
(i)  For each $v\in V_i$, $d_{V_i}(v)\leq 2p^2\sqrt{\varepsilon}|V_i|$.

(ii) For each $v\in V_i$, $d_{V_j}(v)\geq(1-4p^2\sqrt{\varepsilon})|V_j|$.
\end{claim}
{\em Proof.} (i) On contrary, suppose that some vertex $v\in V_1$ satisfies $d_{V_1}(v)>2p^{2}\sqrt{\varepsilon}|V_1|$  without loss of generality. It follows from Lemma \ref{stability} (iv) that $d_{V_i}(v)>2p^{2}\sqrt{\varepsilon}|V_1|$ for $2\le i\le p-1$. Let $U_i$ denote the neighborhood of $v$ in $V_i$. For every $1\leq i\neq j\leq p-1$, it follows from Lemma \ref{stability} (iii) that
there are at most $p^{2}\varepsilon |V_i||V_j|$ non-edges between $U_i$ and $U_j$. Thus
\[
e(U_i, U_j)\geq |U_i||U_j|-p^{2}\varepsilon |V_i||V_j|\geq\left(1-\frac{1}{p^{2}}\right)|U_i||U_j|,
\]
where the second inequality holds since $|V_i|\le 2|V_1|\le\frac{|U_i|}{p^{2}\sqrt{\varepsilon}}$  from Lemma \ref{stability} (ii). If we pick a random vertex from $U_i$ for each $1\leq i\leq p-1$, then they form a copy of $K_{p-1}$ with probability at least $1-{p\choose 2}/p^{2}\geq \frac{1}{2}$. Therefore, the neighborhood of $v$ contains at least
\[
\frac{1}{2}\prod_{i=1}^{p-1} |U_i|\geq \frac{1}{2}\left(2p^{2}\sqrt{\varepsilon}|V_1|\right)^{p-1}\geq  \frac{1}{2}\left(2p^{2}\sqrt{\varepsilon}\right)^{p-1} \left(\frac{1}{p-1}-\sqrt{2\varepsilon}\right)^{p-1} n^{p-1}:=\eta'(p,\varepsilon)n^{p-1}
\]
copies of $K_{p-1}$. By Lemma \ref{Erdos}, this implies that if $\delta$ is sufficiently small in terms of $p$, $t$ and $\varepsilon$, then the neighborhood of $v$ contains a copy of $K_{p-1}(a_2,\dots,a_p)$. Thus, $G$ contains a copy of $K_p(1,a_2,\dots,a_p)$, which is a contradiction.

\medskip
(ii) Suppose that some vertex $v\in V_i$ satisfies $d_{V_j}(v)<(1-4p^2\sqrt{\varepsilon})|V_j|$ for some $j\neq i$. From Claim \ref{ext-fhw} (i), $v$ has at least $(1-2p^{2}\sqrt{\varepsilon})|V_i|-1$ non-neighbors in $V_i$. In total, the number of non-neighbors of $v$ is at least
\[
4p^2\sqrt{\varepsilon}|V_j|+(1-2p^2 \sqrt{\varepsilon})|V_i|-1>\left(1+2p^2\sqrt{\varepsilon}\right)\left(\frac{n}{p-1}-\sqrt{2\varepsilon}n\right)-1>\left(\frac{1}{p-1}+\varepsilon\right)n,
\]
provided $\varepsilon$ is sufficiently small in terms of $p$ and $t$. This contradicts the assumption that the minimum degree of $G$ is at least $(1-\frac{1}{p-1}-\gamma)n$.\hfill$\Box$

\medskip
Suppose the assertion does not hold, then we may assume that there exists some vertex $v\in V_1$ such that $v$ has $a_2$ neighbors in $V_1$ without loss of generality. By Claim \ref{ext-fhw} (ii), we obtain that $v$ and these $a_2$ neighbors have at least $(1-4p^2\sqrt{\varepsilon}(a_2+1))|V_2|>a_3$ common neighbors in $V_3$. Then we can inductively apply Claim \ref{ext-fhw} (ii) to obtain a $K_p(1,a_2,\dots,a_{p})$ in graph $G$, which leads to a contradiction. The assertion follows.\hfill$\Box$

\medskip
The following result is the greedy embedding lemma that we need.
\begin{lemma}[Fox, He and Wigderson \cite{fox}]\label{embed}
Let $k$, $r$, $s$, $t$ be positive integers with $s\leq t$ and $2k\leq t$, and let $H$ be any graph. Let $G$ be an $H$-free graph with $N\geq {t\choose s}^r \frac{t}{2ks}r(H,K_s)$ vertices which contains $K_r(t)$ as an induced subgraph, with parts $V_1,\ldots,V_r$. If $\overline{G}$ does not contain a book $B_{k,n}$ with $n\geq (1-\frac{4ks}{t})\frac{N}{r}$ vertices, then $G$ contains an induced copy of $K_{r+1}(s)$ with parts $W_0,\ldots,W_r$, where $W_i\subseteq V_i$ for every $1\leq i\leq r$.
\end{lemma}

With slightly more careful calculation, we extend a result of \cite[Lemma 2.3]{fox} as follows.
\begin{lemma}\label{fina-embed}
Let $k$, $p$, $x$ be positive integers, and let $z=(4kp)^{p} x$. For any sufficiently small $\mu>0$, there exists an integer $n_0$, such that for all $n\ge n_0$, if $G$ is $H$-free on $N\ge(p-1-\mu)n$ vertices and $\overline{G}$ is $B_{k,n}$-free, and $S\subseteq V(G)$ with $|S|\geq p^{z}\cdot r(H,K_z)$, then $G$ contains $K_{p-1}(x)$ as an induced subgraph, one part of which is a subset of $S$.
\end{lemma}
\noindent{\bf Proof.} We assume $p\ge2$ as it is trivial when $p=1$.
Let $\varepsilon=\frac{1}{4kp}$, $t_{p-1}=x$ and $t_r=\frac{t_{r+1}}{\varepsilon}$ for $r\in[p-2]$. Clearly,
$t_1=t_{p-1}/\varepsilon^{p-2}<z.$
We prove the assertion by induction on $r$ for $r\in [p-1]$.
For $r=1$, we have $|S|>r(H,K_{t_1})$, so $G$ contains an independent set of order $t_1$, that is, $G[S]$ contains $K_r(t_r)$ with $r=1$ as an induced subgraph.
Now suppose $G$ contains $K_r(t_r)$ as an induced subgraph for $r\in[p-2]$, with the first part a subset of $S$. We apply Lemma \ref{embed} with $s=t_{r+1}$ and $t=t_r$ to obtain that
\begin{align*}
{t_r\choose t_{r+1}}^r \frac{t_r}{2kt_{r+1}} r(H,K_{t_{r+1}})\leq \left(\frac{e}{\varepsilon}\right)^{\varepsilon r t_r} 2p\cdot r(H,K_{t_{r+1}})\leq p^{z}\cdot r(H,K_{z})\leq |S|.
\end{align*}
So either $\overline{G}$ contains a $k$-book with at least $(1-4k\varepsilon)\frac{N}{r}\ge(1-\frac1p-\mu)\frac{N}{p-2}\geq n$ vertices for large $n$ and sufficiently small $\mu$, which leads to a contradiction, or $G$ contains an induced $K_{r+1}(t_{r+1})$ whose first $r$ parts are subsets of the $r$ parts of the $K_r(t_r)$. In particular, the first part of this induced $K_{r+1}(t_{r+1})$ is a subset of $S$.
Therefore, the assertion follows inductively.
\hfill$\Box$

\medskip
The following extends a result of  \cite[Lemma 5.1]{fox}.
\begin{lemma}\label{degree-limit}
For every $k, p, t\geq 2$ and for any $\alpha>0$ and sufficiently small $\mu>0$, there exists $\delta>0$ such that the following holds for all large $n$. Let $1\leq a_1\leq \dots \leq a_{p-1}\leq t$ and $a_p\leq\delta n$ be positive integers. If $G$ is a graph on $N\ge (p-1-\mu)n$ vertices such that $G$ is $K_p(a_1,\dots,a_{p})$-free and $\overline{G}$ is $B_{k,n}$-free, then all but at most $\alpha N$ vertices of $G$ have degree at least $(1-\frac{1}{p-1}-\alpha)N$.
\end{lemma}
\noindent{\bf Sketch of the proof.} The idea of the proof is similar to that of \cite[Lemma 5.1]{fox}.
Let $H=K_p(a_1,\dots,a_{p})$. Let $S\subset V(G)$ be the set of vertices of degree less than $(1-\frac{1}{p-1}-\alpha)N$. Let $\varepsilon=\varepsilon(k, \alpha) > 0$, $\delta = \delta(k, p, t)>0$ and $x =\max\{2k,\lceil\frac{t}{\varepsilon}\rceil\}$, and let $z=(4kp)^{p} x$.
We aim to show that $|S|\leq \alpha N$. On contrary, suppose that $|S|\geq \alpha N$. Since $|S|\geq \alpha N\geq p^{z}\cdot r(H,K_z)$ as long as $\delta$ is small enough compared to $\alpha$, by Lemma \ref{fina-embed} we have that there is an induced copy of $K_{p-1}(x)$ in $G$ whose parts are $V_1,\ldots, V_{p-1}$ with $V_1\subseteq S$. Thus, all the vertices of $V_1$ have degree less than $d$. Partition the vertices of $G$ into $p$ parts $U_0,\ldots, U_{p-1}$, where for $i\geq 1$, each vertex in $U_i$ has at most $\varepsilon x$ neighbors in $V_i$, and $U_0$ consists of all vertices with more than $\varepsilon x$ neighbors in each $V_i$.

For the cardinality of $U_i$, $0\leq i\leq p-1$, we have
(i) $|U_0|<(\frac {e}{\varepsilon})^{tp}\cdot \delta n$, (ii) $|U_1|<\frac{n-k}{1-2k\varepsilon}-(\frac{\alpha}{10})^k N$, (iii) $|U_i|<\frac{n-k}{1-2k\varepsilon}$, $i\geq 2$. Indeed, if (i) is false, by probabilistic method, we randomly and uniformly take $a_i$ vertices from $V_i$, $1\leq i\leq p-1$, and take $a_p$ vertices from $U_0$, we can find a copy of $K_p(a_1,\dots,a_{p})$ in $G$. If (ii) or (iii) is false, by probabilistic method, we randomly take $k$ vertices from $U_i$, $1\leq i\leq p-1$, by computation, we know that these $k$ vertices have at least $n-k$ common non-neighbors in $U_i$ and $\overline{U_i}$. So we can find a copy of $B_{k,n}$ in $\overline{G}$. Notice that the difference of (ii) and (iii) is that whether $n-k$ common non-neighbors can be found in $U_i$ or not.
Adding these $U_i$ together, we obtain that the number of vertices in $G$ is less than $N$. This is a contradiction and we are done.
\hfill$\Box$

\section{Proof of Theorem \ref{main}}\label{pf-1}



If $a_1=1$ and $a_2|(n-1-k)$, then the lower bound follows from (\ref{low-bou}), so we will focus on the upper bound in the following.
For every $k, p, t\geq 2$ and for any sufficiently small $\alpha>0$, let $n$ be a sufficiently large integer.
Denote $N=(p-1)(n-1+k(a_2-1))+1$.
For any $K_p(a_1,\dots,a_{p})$-free graph $G$ on $N$ vertices, we aim to show that its complement $\overline{G}$ contains a copy of $B_{k,n}$.

On contrary, suppose that $\overline{G}$ contains no copy of $B_{k,n}$. Lemma \ref{degree-limit} shows that at most $\alpha N$ vertices of $G$ have degree at most $d:=\lceil(1-\frac{1}{p-1}-\alpha)N\rceil.$
Let $T$ be the set of vertices of degree greater than $d$, then the induced subgraph $G[T]$ has at least $(1-\alpha)N$ vertices and thus minimum degree at least $(1-\frac{1}{p-1}-\alpha)N-\alpha N\ge(1-\frac{1}{p-1}-2\alpha)|T|$. Applying Lemma \ref{AES Th} to the graph $G[T]$  which is $K_p(a_1,\dots,a_{p})$-free, since $\alpha$ is sufficiently small in terms of $p$ and $t$, we have that there is a partition $T=\sqcup_{1\le i\le p-1} T_i$, such that each vertex of $T_i$ has at most $a_2-1$ neighbors in $T_i$, for $i=1, 2, \dots, p-1$.

\begin{claim}\label{T-i}
Let $\xi=2p^2 \alpha$. For every $1\leq i\neq j\leq p-1$, we have that

\medskip
(i) \ $ \big||T_i|-\frac{N}{p-1}\big|\le\xi N$.

\smallskip
(ii) \ $d_{T_j}(v)\geq (1-\xi)|T_j|$ for every $v\in T_i$.
\end{claim}

\noindent
{\bf Proof.} (i) Since each vertex of $T_i$ has at most $a_2-1$ neighbors in $T_i$, every vertex in $T_i$ has degree at most $N-|T_i|+a_2-1$. Since every vertex in $T_i$ has degree greater than $d$, this implies that
$|T_i|\leq N-d+a_2-2\le(\frac{1}{p-1}+\alpha)N+a_2-2\leq (\frac{1}{p-1}+2\alpha)N.$
Moreover, since $T_1,\dots,T_{p-1}$ partition $T$, we have
 $|T_i|\geq(1-\alpha)N-(p-2)(\frac{1}{p-1}+2\alpha)N\ge (\frac{1}{p-1}-2p\alpha)N.$

(ii) Note that each vertex of $T_i$ has at most $a_2-1$ neighbors in $T_i$ and $|T_i|\ge (\frac{1}{p-1}-2p\alpha)N$. On contrary, suppose that there is some $v\in T_i$, $d_{T_j}(v)< (1-\xi)|T_j|$. Then
\begin{align*}
d_T(v)< |T|-(|T_i|-a_2)- \xi|T_j|\le\left(1-(1+\xi)\left(\frac{1}{p-1}-2p\alpha\right)\right)|T|+a_2,
\end{align*}
which is at most $(1-\frac{1}{p-1}-2\alpha)|T|$ by noting $\xi=2p^2\alpha$, contradicting to the fact that $G[T]$ has minimum degree at least $(1-\frac{1}{p-1}-2\alpha)|T|$.
\hfill$\Box$

\medskip
Let $S=V(G)\setminus T$, i.e., the set of vertices in $G$ with degree at most $d$. Clearly, $|S|\leq \alpha N$.

\begin{claim}\label{S1-S2}
Let $\zeta=2p^3 t\alpha$. For every $v\in S$, we have that either $d_{T_i}(v)\leq a_2-1$ for some $T_i$, or else there exist two distinct sets $T_k$ and $T_\ell$ such that $d_{T_i}(v)<\zeta|T_i|$ for $i=k,\ell$.
\end{claim}

\noindent
{\bf Proof.} On contrary, suppose that the assertion is false for some $v\in S$. Thus, $d_{T_i}(v)\geq \zeta|T_i|$ for all but at most one choice of $i\in[p-1]$, and $d_{T_i}(v)\geq a_2$ for each $T_i$. By relabeling the parts, we may assume that $d_{T_i}(v)\geq \zeta|T_i|$ for all $i\in[p-2]$. We obtain an $a_2$-set $M$ by picking $a_2$ neighbors of $v$ in $V_{p-1}$. By Claim \ref{T-i} (ii), we see that $v$ and all vertices of $M$ have at least $(\zeta-a_2\xi)|T_1|> \xi |T_1|>a_3$ common neighbors in $T_1$, for $N$ sufficiently large. Continuing in this way, we can get a copy of $K_p(a_1,\dots,a_{p})$ in $G$ in the end, which leads to a contradiction. \hfill$\Box$

\medskip
From the above claim, we first partition $S$ into two subsets $S_1$ and $S_2$, where $S_1$ consists of all vertices in $S$ that each has at most $a_2-1$ neighbors in some $T_i$, and $S_2$ consists of the remaining vertices $v$, namely those satisfying $d_{T_i}(v)<\zeta|T_i|$ for at least two choices of $1\leq i\leq p-1$. Second, we partition $S_1$ into $p-1$ subsets $S_{11},\dots,S_{1,p-1}$ such that a vertex $v\in S_{1}$ is put into $S_{1i}$ for $1\leq i\leq p-1$ if and only if $i$ is the smallest index such that $d_{T_i}(v)\leq a_2-1$. Finally, we partition $S_2$ into $p-1$ subsets $S_{21},\dots,S_{2,p-1}$ such that a vertex $v\in S_{2}$ is put into $S_{2i}$ for $1\leq i\leq p-1$ if and only if $i$ is the smallest index satisfying $d_{T_i}(v)<\zeta|T_i|$.
In the following, we consider two cases according to the size of $S_2$.

In the first case, we assume that $|S_2|< 2p$. For this case, we define $U_i=T_i\cup S_{1i}\cup S_{2i}$, so we obtain that $V=\sqcup_{i=1}^{p-1} U_i$. Therefore, there exists some $h\in [p-1]$ such that $$|U_h|\geq \left\lceil\frac{|V|}{p-1}\right\rceil\ge n+k(a_2-1).$$

By Lemma \ref{AES Th}, any vertex $x\in {T_h}$ satisfies $d_{T_h}(x)\le a_2-1$. Therefore, the Tur\'{a}n's bound (Lemma \ref{turan}) implies that there exists an independent set $X\subseteq T_h$ with $|X|\ge \frac{|T_h|}{a_2}$.
Note that $|S_{2h}|\leq |S_2|<2p$ and any vertex $v\in S_{1h}$ satisfies $d_{T_h}(v)< a_2$, so we can find an independent set $X'\subset X$ such that any vertex of $X'$ is non-adjacent to all vertices of $S_{1h}\cup S_{2h}$ with
\begin{align*}
|X'|\ge |X|-a_2|S_{1h}|-\zeta |S_{2h}| |T_h|\ge \frac{|T_h|}{a_2}-t\alpha N-2p\zeta |T_h|.
\end{align*}
Note from Claim \ref{T-i} (i) that $|T_h|\geq (\frac{1}{p-1}-\xi)N$, we obtain that
\begin{align*}
|X'|\ge\left(\frac{1}{p-1}-\xi\right)\left(\frac{1}{a_2}-2p\zeta\right)N-t\alpha N \geq k
\end{align*}
by noting $\xi$, $\zeta$ and $\alpha$ are sufficiently small in terms of $p$, $a_2$ and $t$, and $n$ is large.

Now, applying Lemma \ref{AES Th} again that any vertex $x\in X$ satisfies $d_{T_h}(x)\le a_2-1$, we can find an independent set of order $k$ in $X'$ such that it is non-adjacent to all but at most $k(a_2-1)$ vertices of $T_h$. Consequently, these $k$ vertices have at least
\[
|T_h\cup S_{1h}\cup S_{2h}|-k-k(a_2-1)=|U_h|-k-k(a_2-1)\ge n-k
\]
common non-neighbors, i.e., there exists a $B_{k,n}$ in $\overline{G}$.

In the second case, we may assume $|S_{2}|\ge 2p$. Note that $S_2=\sqcup_{1\le i\le p-1}S_{2i}$, where $S_{2i}$ consists of all vertices $v\in S_2$ such that $i$ is the smallest index satisfying $d_{T_i}(v)<\zeta|T_i|$. Now we partition $S_2$ in another way. Indeed, we partition $S_2$ into $p-1$ subsets $S_{21}',\dots,S_{2,p-1}'$ such that a vertex $v\in S_{2}$ is put into $S_{2i}'$ for $1\leq i\leq p-1$ if and only if $i$ is the biggest index satisfying $d_{T_i}(v)<\zeta|T_i|$.

Since $S_2$ consists of vertices $v\in S$ satisfying $d_{T_i}(v)<\zeta|T_i|$ for at least two choices of $1\leq i\leq p-1$, we have that $S_{2i}\cap S_{2i}'=\emptyset$.  Clearly, $|S_{2i}|+|S_{2i}'|\leq |S_2|$.
We define $$W_i=T_i\cup S_{1i}\cup S_{2i}\cup S_{2i}'.$$
Note that $(T_i\cup S_{1i}\cup S_{2i})\cap(T_j\cup S_{1j}\cup S_{2j})=\emptyset$ for distinct $i$ and $j$, and $S_{2i}' \cap(T_i\cup S_{1i}\cup S_{2i})=\emptyset$, so we obtain that
\begin{align*}
\sum_{i=1}^{p-1} |W_i|=\sum_{i=1}^{p-1} |T_i\cup S_{1i}\cup S_{2i}\cup S_{2i}'|
=\sum_{i=1}^{p-1}\big (|T_i\cup S_{1i}\cup S_{2i}|+|S_{2i}'|\big)=|V|+|S_2|.
\end{align*}
It follows by the pigeon-hole principle that there exists some $h\in [p-1]$ such that
\begin{align}\label{Wh}
|W_h|=|T_h\cup S_{1h}\cup S_{2h}\cup S_{2h}'|\geq \frac{|V|+|S_2|}{p-1}\ge n-1+\frac{|S_2|}{p-1}+k(a_2-1).
\end{align}
Recall that any vertex $x\in {T_h}$ has at most $a_2-1$ neighbors in $T_h$. Therefore, by Lemma \ref{turan}, we see that there exists an independent set $X\subseteq T_h$ with $|X|\ge \frac{|T_h|}{a_2}$. We randomly and uniformly pick a $k$-set $K\subset X$, which is an independent set. We will count the number of common non-neighbors of vertices of $K$.

First, all vertices of $K$ have at least $|T_h|-ka_2$ common non-neighbors in $T_h$ since any vertex $v\in {T_h}$ has at most $a_2-1$ neighbors in $T_h$. Next, since any vertex $v\in {S_{2h}\cup S_{2h}'}$ satisfies $d_{X}(v)\le d_{T_h}(v)< \zeta |T_h|<a_2\zeta|X|$, it follows that the probability of any vertex in ${S_{2h}\cup S_{2h}'}$ can be chosen as a common non-neighbor of $K$ is at least $1-ka_2\zeta$. Finally, if $v\in S_{1h}$, the probability of any vertex in $S_{1h}$ can be chosen as a common non-neighbors of $K$ is at least $1-k\frac{a_2-1}{|X|}$, since any vertex $v\in {S_{1h}}$ satisfies $d_{X}(v)\le d_{T_h}(v)\le a_2-1$.  By the linear of expectation, we conclude that there exists a $k$-set $K_0$ in $X$ such that the number $\lambda$ of common non-neighbors of vertices of $K_0$ in $T_h\cup S_{1h}\cup S_{2h}\cup S_{2h}'$ satisfies
\begin{align*}
\lambda\ge\big(|T_h|-ka_2\big)+\big(|S_{2h}|+|S_{2h}'|\big)\big(1-ka_2\zeta\big)+|S_{1h}|\left(1-k\frac{a_2-1}{|X|}\right).
\end{align*}
Note that  $\frac{|S_{1h}|}{|X|}\le\frac{\alpha N }{|T_h|/a_2}\le\frac{\alpha  a_2}{1/p}=pa_2\alpha$,
so we have that
\begin{align*}
\lambda \geq\big(|T_h|+|S_{2h}|+|S_{2h}'|+|S_{1h}|\big)-ka_2-ka_2\zeta\big(|S_{2h}|+|S_{2h}'|\big)-kpa_2(a_2-1)\alpha.
\end{align*}
Therefore, by inequality (\ref{Wh}) and $|S_{2h}|+|S_{2h}'|\leq |S_2|$, we obtain that
\begin{align*}
\lambda \geq& n-1+ \frac{|S_2|}{p-1}+k(a_2-1)-ka_2-ka_2\zeta |S_2|-kpa_2(a_2-1)\alpha
\\ \geq & n-1+ \left(\frac{1}{p-1}-ka_2\zeta\right) |S_2|-k-kpa_2(a_2-1)\alpha
\\\ge & n-k
\end{align*}
where the last inequality follows by noting $|S_2|\ge 2p$, and $\zeta$ and $\alpha$ are sufficiently small. Thus we can get a $B_{k,n}$ in $\overline{G}$.

\medskip
Now, combining the above two cases, we can get the final contradiction. This completes the proof of Theorem \ref{main}.
\hfill$\Box$

\section{Proof of Theorem \ref{main-2}}\label{pf-2}

 Erd\H{o}s and Stone \cite{e-stone} proved that every graph with $n$ vertices and at least $(1-\frac{1}{p}+\varepsilon)\frac{n^2}2$ edges contains a $K_{p+1}(t)$, where, for $p$ and $\varepsilon$ fixed, $t=t(p,\varepsilon,n)$ tends to infinity with $n$. In \cite{be-73,bes,bk,c-sz,i}, the function $t$ has been widely studied, and we may take $t=\Omega(\ln n)$ for fixed $p$ and $\varepsilon$. We shall apply  the following  result by Nikiforov  \cite{ni} which concludes with the same assertion from weaker premises.
\begin{lemma}[Nikiforov \cite{ni}]\label{niki}
For every $p\geq 2$, $\eta>0$ there exists $\delta>0$ such that if $G$ has $n$ vertices and at least $\eta n^p$ copies of $K_p$, then $K_p(t)\subset G$ for $t=\lceil \delta\ln n\rceil$.
\end{lemma}

We also need the following result by Bollob\'{a}s and Nikiforov \cite{bn}, in which the minimum degree condition can be weakened to the corresponding average degree condition.
\begin{lemma}[Bollob\'{a}s and Nikiforov \cite{bn}]\label{bn}
Let $p\geq 1$ and $\alpha\geq 0$. If $G$ has $n$ vertices and minimum degree $\delta(G)\geq(1-\frac{1}{p}+\alpha)n$, then $G$ contains at least $\alpha\frac{p^2}{p+1}(\frac{n}{p})^{p+1}$ copies of $K_{p+1}$.
\end{lemma}

We need a  property for a graph that is $K_p(1,a,b,\dots,b)$-free with high minimum degree as follows, where $b\leq \delta\ln n$.

\begin{lemma}\label{AES Th-2}
For every $a\geq 1$, and $p\geq 2$, there exists $\gamma,\delta>0$ such that if $n$ is large enough in terms of $a$ and $p$, $b\leq \delta\ln n$, and $G$ is a $K_p(1,a,b,\dots,b)$-free graph on $n$ vertices with minimum degree at least $(1-\frac{1}{p-1}-\gamma)n$, then there is a partition $V(G)=\sqcup_{1\le i\le p-1} V_i$ such that each vertex of $V_i$ has at most $a-1$ neighbors in $V_i$, for $i=1, 2, \dots, p-1$.
\end{lemma}
\noindent{\bf Proof.} 
For $\varepsilon>0$, where $\varepsilon$ is sufficiently small in terms of $p$ and $a$, let $\gamma$ and $\eta$ be applied in Lemma \ref{stability}. Further, let $\delta>0$ be the parameter guaranteed from Lemma \ref{niki}. It is clear that $G$ is $K_p(\lceil \delta\ln n\rceil)$-free from the assumption, so we have that $G$ has at most $\eta n^p$ copies of $K_p$ by Lemma \ref{niki}. Since $G$ has at most $\eta n^p$ copies of $K_p$ with minimum degree at least $(1-\frac{1}{p-1}-\gamma)n$, it follows from  Lemma \ref{stability} that there is a partition $V(G)=\sqcup_{1\le i\le p-1} V_i$ such that (i)-(iv) of Lemma \ref{stability} hold. By a similar argument as in Lemma \ref{AES Th}, we have the following claim.

\begin{claim}\label{ext-fhw-2}
For distinct $i$ and $j$, $1\leq i,j\leq p-1$, we have

\medskip
(i)  For each $v\in V_i$, $d_{V_i}(v)\leq 2p^2\sqrt{\varepsilon}|V_i|$.

(ii) For each $v\in V_i$, $d_{V_j}(v)\geq(1-4p^2\sqrt{\varepsilon})|V_j|$.
\end{claim}

Suppose the assertion does not hold, then we may assume that there exists some vertex $v\in V_1$ such that $d_{V_1}(v)\geq a$ without loss of generality. Let us take a subset $A$ with $|A|=a$ of the neighborhood of $v$ in $V_1$, and let $U_j$ denote the common neighborhood of $v$ and $A$ in $V_j$, where $2\leq j\leq p-1$. By Claim \ref{ext-fhw-2} (ii), we obtain that $|U_j|\geq (1-4p^2 \sqrt \varepsilon(a+1))|V_j|\geq \frac{1}{2}|V_j|$.
For every $2\leq i\neq j\leq p-1$, we have by Lemma \ref{stability} (iii) that
\[
e(U_i, U_j)\geq |U_i||U_j|-p^{2}\varepsilon |V_i||V_j|\geq \left(1-4p^{2}\varepsilon\right)|U_i||U_j|\geq\left(1-\frac{1}{(p-2)^{2}}\right)|U_i||U_j|,
\]
where the second inequality uses $|U_i|\geq \frac{1}{2}|V_i|$. By the union bound, if we pick a random vertex from $U_j$ for each $2\leq j\leq p-1$, then they span a copy of $K_{p-2}$ with probability at least $1-{p-2\choose 2}/(p-2)^{2}\geq \frac{1}{2}$. Therefore, the common neighborhood of $v$ and $A$ contains at least
\[
\frac{1}{2}\prod_{j=2}^{p-1} |U_j|\geq \frac{1}{2^{p-1}}\prod_{j=2}^{p-1} |V_j|\geq \frac{1}{2^{p-1}p^{p-2}} n^{p-2}:=\eta'n^{p-2}
\]
copies of $K_{p-2}$, for $\eta'$ depending on $p$. By Lemma \ref{niki}, we have that the common neighborhood of $v$ and $A$ contains a copy of $K_{p-2}(b)$, so $G$ contains a copy of $K_p(1,a,b,\dots,b)$, which leads to a contradiction. Thus the assertion follows.\hfill$\Box$


\begin{lemma}\label{degree-limit-2}
For every $k, a\geq 1$, $p\geq 2$ and for any $\alpha>0$ and sufficiently small $\mu>0$, there exists $\delta>0$ such that the following holds for all large $n$. If $G$ is a graph on $N\ge (p-1-\mu)n$ vertices such that $G$ is $K_p(1,a,b,\dots,b)$-free with $b\leq \delta\ln n$ and $\overline{G}$ is $B_{k,n}$-free, then all but at most $\alpha N$ vertices of $G$ have degree at least $(1-\frac{1}{p-1}-\alpha)N$.
\end{lemma}
\noindent{\bf Sketch of the proof.} The proof is similar to that of Lemma \ref{degree-limit}.
Let $H=K_p(1,a,b,\dots,b)$. Let $S\subset V(G)$ be the set of vertices of degree less than $(1-\frac{1}{p-1}-\alpha)N$. Let $\varepsilon=\varepsilon(k, \alpha) > 0$, $\delta = \delta(k, p, t)>0$ and $x =\max\{\lceil\frac{b}{\varepsilon}\rceil, 2k\}$, and let $z=(4kp)^{p} x$. We aim to show that $|S|\leq \alpha N$. On contrary, suppose that $|S|\geq \alpha N$. Note that $r(K_t,K_t)<4^t$ (see \cite{esz}), so we have
 $$p^{z}\cdot r(H,K_z)\le p^{(4kp)^p x}
 \cdot  r(K_{pb},K_{(4kp)^p x}) \leq \alpha N\le|S|$$ as long as $\delta$ is small enough compared to $\alpha$ and $\varepsilon$. Therefore, it follows by Lemma \ref{fina-embed} that there is an induced copy of $K_{p-1}(x)$ in $G$ whose parts are $V_1,\ldots, V_{p-1}$ with $V_1\subseteq S$. Thus, all the vertices of $V_1$ have degree less than $d$. Partition the vertices of $G$ into $p$ parts $U_0,\ldots, U_{p-1}$, where for $i\geq 1$, each vertex in $U_i$ has at most $\varepsilon x$ neighbors in $V_i$, and $U_0$ consists of all vertices with more than $\varepsilon x$ neighbors in each $V_i$.

Similar to the proof in Lemma \ref{degree-limit}, for $0\leq i\leq p-1$, we have
(i) $|U_0|<(\frac {e}{\varepsilon})^{pb}\cdot b$, (ii) $|U_1|<\frac{n-k}{1-2k\varepsilon}-(\frac{\alpha}{10})^k N$, (iii) $|U_i|<\frac{n-k}{1-2k\varepsilon}$, $i\geq 2$.
Adding these $U_i$ together, we obtain that the number of vertices in $G$ is less than $N$. This is a contradiction and we are done.
\hfill$\Box$

\medskip
Now we are ready to prove Theorem \ref{main-2}.

\medskip
\noindent
{\bf Proof of Theorem \ref{main-2}} \; If $a|(n-1-k)$, then the lower bound follows from (\ref{low-bou}), so we will focus on the upper bound in the following. Note that  the assertion holds from Theorem \ref{main} for $p=2,3$.
For every $k, a\geq 1$, $p\geq 4$ and for any sufficiently small $\alpha>0$, let $n$ be a sufficiently large integer.
Let $N=(p-1)(n-1+k(a-1))+1$.
It suffices to show that for any $K_p(1,a,b,\ldots,b)$-free graph $G$ on $N$ vertices,  its complement $\overline{G}$ contains a copy of $B_{k,n}$.

On contrary, suppose that $\overline{G}$ contains no copy of $B_{k,n}$. Lemma \ref{degree-limit-2} shows that at most $\alpha N$ vertices of $G$ have degree at most $d:=\lceil(1-\frac{1}{p-1}-\alpha)N\rceil.$
Let $T$ be the set of vertices of degree greater than $d$, then the induced subgraph $G[T]$ has at least $(1-\alpha)N$ vertices and thus minimum degree at least $(1-\frac{1}{p-1}-2\alpha)|T|$. Applying Lemma \ref{AES Th-2} to the graph $G[T]$  which is $K_p(1,a,b,\ldots,b)$-free, since $\alpha$ is sufficiently small in terms of $p$ and $a$, we have that there is a partition $T=\sqcup_{1\le i\le p-1} T_i$, such that each vertex of $T_i$ has at most $a-1$ neighbors in $T_i$, for $i=1, 2, \dots, p-1$.

By a similar argument as Claim \ref{T-i} using Lemma \ref{AES Th-2}, we have the following claim.

\begin{claim}\label{T-i-2}
Let $\xi=2p^2 \alpha$. For every $1\leq i\neq j\leq p-1$, we have that

\medskip
(i) \ $ \big||T_i|-\frac{N}{p-1}\big|\le\xi N$.

\smallskip
(ii) \ $d_{T_j}(v)\geq (1-\xi)|T_j|$ for every $v\in T_i$.
\end{claim}

Let $S=V(G)\setminus T$, i.e., the set of vertices in $G$ with degree at most $d$, and $|S|\leq \alpha N$.

\begin{claim}\label{S1-S2}
Let $\zeta=8p^4a\alpha$. For every $v\in S$, we have that either $d_{T_i}(v)\leq a-1$ for some $T_i$, or else  there exist two distinct sets $T_k$ and $T_\ell$ such that $d_{T_i}(v)<\zeta|T_i|$ for $i=k,\ell$.
\end{claim}

\noindent
{\bf Proof.} On contrary, suppose that the assertion is false for some $v\in S$. Thus, $d_{T_i}(v)\geq \zeta|T_i|$ for all but at most one choice of $i\in[p-1]$, and $d_{T_i}(v)\geq a$ for each $T_i$. By relabeling the parts, we may assume that $d_{T_i}(v)\geq \zeta|T_i|$ for all $i\in[p-2]$. Let $A$ denote the neighborhood of $v$ in $T_{p-1}$. Let $W_i$ denote the common neighborhood of $v$ and $A$ in $T_i$. By Claim \ref{T-i-2} (ii), we obtain that $|W_i|\geq (\zeta-a\xi)|T_i|\geq \frac{\zeta}{2}|T_i|$. Let $W=\cup_{i=1}^{p-2} W_i$, and let $w\in W$ be a vertex of minimum degree in $G[W]$, say $w\in W_{i_0}$ for some $i_0\in[p-2]$. For $p\ge4$,  we have that
\begin{align*}
\delta(G[W])\ge&\sum_{1\le j\le p-2, j\neq i_0} d_{W_j}(w)
\geq \sum_{1\le j\le p-2, j\neq i_0} \left(d_{T_j}(w)+|W_j|-|T_j|\right)
          \\\geq &\frac{p-3}{p-2}\left(1-\frac{2\xi}{\zeta}\right)|W|\geq \left(\frac{p-4}{p-3}+\frac{1}{2(p-2)(p-3)}\right)|W|.
\end{align*}
Apply Lemma \ref{bn} to $G[W]$, we obtain that the number of $K_{p-2}$ in $G[W]$ is at least
$$
\frac{1}{2(p-2)(p-3)}\frac{(p-3)^2}{p-2}\left(\frac{|W|}{p-3}\right)^{p-2}\ge
\frac{1}{2(p-2)(p-3)}\frac{(p-3)^2}{p-2}\left(\frac{\zeta(p-2)|T_i|}{2(p-3)}\right)^{p-2},
$$
which is at least $\eta N^{p-2}$ provided $\eta$ is sufficiently small compared with $\zeta$.
By Lemma \ref{niki}, we have that the common neighborhood of $v$ and $A$ contains a copy of $K_{p-2}(b)$. So $G$ contains a copy of $K_{p}(1,a,b,\ldots,b)$, which leads to a contradiction. \hfill$\Box$

\medskip
Now we can complete the proof of Theorem \ref{main-2} as that of Theorem \ref{main}.
\hfill$\Box$

\section{Concluding remarks}\label{clu}

We have known from Section \ref{disp} that Conjecture \ref{conj} does not hold for many cases. However, Conjecture \ref{conj} is still open for $a_1\ge2$. In particular, it would be interesting to compute more exact values for $r(C_4,B_n)$, although we have obtained the value of $r(C_4,B_n)$ for infinitely many $n$, see \cite{llp}. 
Moreover, as pointed out by Fox, He and Wigderson \cite{fox}, it would be very interesting to see how far one can push these ideas; for example, is it possible to completely eliminate the use of the regularity lemma from the proof of \cite[Theorem 2.1]{nr09}?
One could extend Theorem \ref{main-2} to somewhat general cases using the technique we use combining with a stronger result of Nikiforov \cite[Theorem 1]{ni}. However,  new ideas are needed if we want to eliminate the use of the regularity lemma from the proof of \cite[Theorem 2.1]{nr09}.
On the other direction, as a special case, the behavior of $r(K_n,B_n)$ is much different. Li and Rousseau \cite{lr} obtained that $c_1n^3/\ln^2 n\le r(K_n,B_n)\le c_2n^3/\ln n$ for sufficiently large $n$, where $c_1$ and $c_2$ are positive constants.
Sudakov \cite{sud} improved the upper bound by a factor $\sqrt{\ln n}$, and also conjectured that $r(K_n,B_n)=\Theta(n^3/\ln^2 n)$.


\end{spacing}

\end{document}